\documentclass[draft,preprint,12pt,numbers,sort&compress]{elsarticle}
\usepackage{}
\usepackage{mathrsfs}
\usepackage{amssymb}
\usepackage{amsthm}
\usepackage{mathrsfs}
\usepackage[centertags]{amsmath}
\usepackage{amsfonts}

\usepackage{color}

\usepackage{ctex}
\usepackage{CJK}

\input amssym.def
\input amssym.tex

\date{}

\newtheorem{Theorem}{Theorem}[section]

\newtheorem{Lemma}{Lemma}[section]

\newcommand\R{\mbox{\bf R}}

\newcommand\T{\mbox{\bf T}}
\newcommand\SR{\mbox{\scriptsize\bf R}}

\newcommand{\definition}{{\lower .5ex
  \hbox{$\>\>\stackrel{\triangle}{=}\>\>$} }}

\begin{document}

\baselineskip=22pt
\thispagestyle{empty}
\mbox{}
\bigskip

\begin{center}
{\Large \bf Convergence problem of the Kawahara equation on the real line}\\[1ex]

{ Wei Yan
\footnote{Email:011133@htu.edu.cn}$^a$, \qquad  Weimin Wang \footnote{Email: ydn1129@163.com}$^{a}$, \qquad  Xiangqian Yan \footnote{Email: yanxiangqian213@126.com}$^{a}$}\\[1ex]

{$^a$School of Mathematics and Information Science, Henan
Normal University,}\\
{Xinxiang, Henan 453007,   China}\\[1ex]

\end{center}
\noindent{\bf Abstract.}
In this paper, we consider  the
convergence
 problem of  the Kawahara equation
 \begin{eqnarray*}
 &&u_{t}+\alpha\partial_{x}^{5}u+\beta\partial_{x}^{3}u+\partial_{x}(u^{2})=0
 \end{eqnarray*}
on the real line  with rough data. Firstly,  by using Strichartz estimates as well as high-low frequency idea,
 we establish two crucial bilinear estimates, which are just Lemmas 3.1-3.2 in this paper;  we also present the proof of Lemma 3.3
which shows that $s>-\frac{1}{2}$ is necessary for   Lemma 3.2.
Secondly, by using frequency truncated technique and high-low frequency technique,  we show the pointwise convergence of
 the Kawahara equation with rough data in $H^{s}(\R)(s\geq\frac{1}{4})$;   more precisely, we prove
 \begin{eqnarray*}
 &&\lim\limits_{t\rightarrow0}u(x,t)=u(x,0), \qquad a.e. x\in\R,
 \end{eqnarray*}
 where $u(x,t)$ is the solution to the Kawahara equation with initial data $u(x,0).$
 Lastly, we show
 \begin{eqnarray*}
 &&\lim\limits_{t\rightarrow0}\sup\limits_{x\in\SR}|u(x,t)-U(t)u_{0}|=0
 \end{eqnarray*}
 with   rough data in $H^{s}(\R)(s>-\frac{1}{2})$.

 \medskip


\medskip
\noindent {\bf Corresponding Author:}Wei Yan

\medskip
\noindent {\bf Email Address:}011133@htu.edu.cn


\leftskip 0 true cm \rightskip 0 true cm

\newpage

\baselineskip=20pt

\bigskip
\bigskip

{\large\bf 1. Introduction}
\medskip

\setcounter{Theorem}{0} \setcounter{Lemma}{0}\setcounter{Definition}{0}\setcounter{Proposition}{1}

\setcounter{section}{1}
In this paper,  we  investigate  the Cauchy problem for the Kawahara equation
\begin{eqnarray}
&&u_{t}+\alpha\partial_{x}^{5}u+\beta\partial_{x}^{3}u+\partial_{x}(u^{2})=0, x\in \R,\label{1.01}\\
&&u(x,0)=u_{0}(x)   \label{1.02},
\end{eqnarray}
where $\alpha\neq0$, $\beta$ are real numbers.

The Kawahara equation arises in the study of the water waves with surface tension,
 in which the Bond number takes on the critical value, where the Bond number represents a
dimensionless magnitude of surface tension in the shallow water regime \cite{BS,K,KO}.
Some authors have  studied  the Cauchy problem for the Kawahara equation on the real line
\cite{BGZ, CDT,CT,CK,H2005,JH2009,K2013,K2011,WCD,YL,YL2011}
and \cite{K2012} on the torus and \cite{CK} on the half-line.  By using the Fourier restriction norm method introduced in \cite{Bourgain1993} and developed in
 \cite{KPV-Duke,KPV1996},
Chen et al. \cite{CLMW} and Jia and Huo \cite{JH2009} independently
proved that the Cauchy problem for (\ref{1.01}) is locally well-posed in $H^{s}(\R)(s>-\frac{7}{4})$. Chen and Guo \cite{CG} proved that the Cauchy problem for (\ref{1.01})
      is globally well-posed in $H^{s}(s\geq -\frac{7}{4})$ with the aid of I-method introduced in  \cite{CKSTT}.
Kato \cite{K2011} proved that the Cauchy problem for the
  Kawahara equation is locally well-posed  in $H^{s}(\R)(s\geq-2)$ with the aid of some resolution
spaces and ill-posed in $H^{s}(\R)(s<-2)$ by using the argument of \cite{BT}. Kato \cite{K2013} proved that the Cauchy problem for the
 Kawahara equation is globally well-posed  in $H^{s}(\R)(s>-\frac{38}{21})$ with the aid of I-method.

 Now we recall the research history of pointwise convergence problem of some linear dispersive equations.
 Carleson \cite{Carleson} firstly investigated    the pointwise convergence  problem of  one dimensional Schr\"odinger equation in $H^{s}(\R)(s\geq 1/4)$.
  Some  authors have investigated  the pointwise
  convergence problem of the  Schr\"odinger equation  in   dimensions $n\geq2$  \cite{Bourgain1992,Bourgain2013,Bourgain2016,CLV,Cowling,DK,DG,Du,DGZ,DZ,GS2008,Lee,
       LR2015,LR2017,LR2019, MVV,RVV,RV,Shao,S,Vega}  and    established an improved maximal inequality for $2D$
   fractional order Schr\"odinger operators \cite{MYZ2015} and    established  the maximal estimates for Schr\"odinger
    equation with inverse-square potential \cite{MZZ2015}.
 Dahlberg and Kenig \cite{DK} showed
that the pointwise convergence problem of the
Schr\"odinger equation is invalid in $H^{s}(\R^{n})(s<\frac{1}{4}).$
Bourgain \cite{Bourgain2016} presented counterexamples showing that
 when $s<\frac{n}{2(n+1)}(n\geq2)$,   the pointwise
 convergence problem of $n$ dimensional Schr\"odinger  equation does not hold in $H^{s}(\R^{n})$.
 Recently, Du et al. \cite{DGL} proved that the pointwise convergence problem of
  two dimensional
  Schr\"odinger
 equation is valid in $H^{s}(\R^{2})$ with $s>\frac{1}{3}$.
 Du and Zhang \cite{DZ} showed that the
  pointwise convergence problem of $n$ dimensional Schr\"odinger
 equation   is valid  for data in $H^{s}(\R^{n})(s>\frac{n}{2(n+1)},n\geq3).$

 Compaan \cite{Compaan} studied the smooth  property and dispersive blow-up of semilinear Schr\"odinger equation.
   Compaan et al. \cite{CLS}
  showed the convergence  problem of
  the  nonlinear Schr\"odinger flows with
   rough data and random data.
   Linares and  Ramos \cite{LR,LR2021}
  showed the pointwise convergence results for the flow of the generalized
   Zakharov-Kuznetsov equation.

 In this paper, motivated by \cite{Compaan, CLS,LR,LR2021},  we  show
  $\lim\limits_{t\rightarrow0}u(x,t)=u(x,0), \qquad a.e. x\in\R, $
  where $u(x,t)$ is the solution to  the Kawahara equation
 with rough  data in $H^{s}(\R)(s\geq\frac{1}{4})$; we also show
  \begin{eqnarray*}
  &&\lim\limits_{t\rightarrow0}\sup\limits_{x\in\SR}|u(x,t)-U(t)u_{0}|=0
  \end{eqnarray*}
  with   rough data in $H^{s}(\R)(s>-\frac{1}{2})$.

We present some notations before stating the main results. $|E|$ denotes by the Lebesgue measure
 of set $E$.
 $A\sim B$ means that there exists $C>0$ such that $\frac{1}{C}|A|\leq |B|\leq C|A|.$
  We define
 a smooth jump function $\eta(\xi)$ such that $\eta(\xi)=1$ for $|\xi|\leq1$ and $\eta(\xi)=0$
  for $|\xi|>2$ and $\phi(\xi)=\alpha\xi^{5}-\beta\xi^{3}$.
 We define
\begin{eqnarray*}
&&a={\rm max}\left\{1,\left(2\left|\frac{3\beta}{5\alpha}\right|\right)^{\frac{1}{2}}\right\},\\
&&\mathscr{F}_{x}f(\xi)=\frac{1}{\sqrt{2\pi}}\int_{\SR}e^{-ix \xi}f(x)dx,\\
&&\mathscr{F}_{x}^{-1}f(\xi)=\frac{1}{\sqrt{2\pi}}\int_{\SR}e^{ix \xi}f(x)dx,\\
&&  \mathscr{F}f(\xi,\tau)=\frac{1}{2\pi}\int_{\SR^{2}}e^{-ix \xi-it\tau}f(x,t)dxdt,\\
 &&  \mathscr{F}^{-1}f(\xi,\tau)=\frac{1}{2\pi}\int_{\SR^{2}}e^{ix \xi+it\tau}f(x,t)dxdt,\\
&&U(t)f=\frac{1}{\sqrt{2\pi}}\int_{\SR} e^{ix\xi-it\phi(\xi)}\mathscr{F}_{x}u_{0}(\xi)d\xi,\\
&&P^{a}f=\frac{1}{2\pi}\int_{|\xi|\geq a}e^{ix\xi}\mathscr{F}_{x}f(\xi)d\xi,\\
&&P_{a}f=\frac{1}{2\pi}\int_{|\xi|\leq a}e^{ix\xi}\mathscr{F}_{x}f(\xi)d\xi.
\end{eqnarray*}
The space $H^{s}(\R)$ is the completion of the Schwartz function space on $\R$ with respect to the norm
$\|f \|_{H^{s}(\SR)}=
\|\langle\xi\rangle^{s}\mathscr{F}_{x}f\|_{L_{\xi}^{2}(\SR)}$, where
 $\langle\xi\rangle^{s}=(1+|\xi|^{2})^{\frac{s}{2}}$. The  space $X_{s,b}(\R^{2})$  is defined to
  be the completion of the Schwartz function space on $\R^{2}$ with respect
 to the norm
\begin{eqnarray*}
\left\|u\right\|_{X_{s,b}(\SR^{2})}=\left\|\langle\xi\rangle^{s}\langle\sigma\rangle^{b}\mathscr{F}u(\xi,\tau)
\right\|_{L_{\xi\tau}^{2}}.
\end{eqnarray*}
Here, $\sigma=\tau+\phi(\xi)$.

The main results are as follows:

\begin{Theorem}\label{Theorem1} (Bilinear estimate related to Kawahara equation)
For $s\geq-\frac{7}{4}+4\epsilon$, $b^{\prime}=-\frac{1}{2}+2\epsilon$
 and $b=\frac{1}{2}+\epsilon$. Then, we have
\begin{eqnarray}
&&\left\|\partial_{x}\left(u_{1}u_{2}\right)\right\|_{X_{s,b^{\prime}}}\leq
 C\left\|u_{1}\right\|_{X_{s,b}}\left\|u_{2}\right\|_{X_{s,b}}.\label{1.03}
\end{eqnarray}
\end{Theorem}
 \noindent{\bf Remark 1.} Theorem 1 is the new proof of Corollary 3.3 of \cite{JH2009}.
  Jia and Huo used the Cauchy-Schwarz inequality and Strichartz
 estimates to establish Corollary 3.3 of \cite{JH2009}.  In this paper, we only use the Strichartz estimates to establish Theorem 1.1. Moreover,
 Lemma 2.6 plays an important role in establishing Theorem 1.1.

\begin{Theorem}\label{Theorem2}  (Bilinear estimate related to Kawahara equation). Let $s_{1}=\frac{1}{2}+2\epsilon$, $b^{\prime}=-\frac{1}{2}+2\epsilon$,
  $b=\frac{1}{2}+\frac{\epsilon}{2}$ and  $s_{2}\geq-\frac{1}{2}+\epsilon$.  Then, we have
\begin{eqnarray}
&&\|\partial_{x}(u_{1}u_{2})\|_{X_{s_{1},b^{\prime}}}\leq C
\|u_{1}\|_{X_{s_{2},b}}\|u_{2}\|_{X_{s_{2},b}}.\label{1.04}
\end{eqnarray}
\end{Theorem}
\noindent{\bf Remark 2.} From the proof process of Theorem 1.5,    Theorem 1.2  plays an import part in establishing Theorem 1.5.

\begin{Theorem} \label{Theorem3} (Necessity of $s>-\frac{1}{2}$  in Theorem 2)
Let $s\leq-\frac{1}{2}$, $s_{1}=\frac{1}{2}+2\epsilon$, $b^{\prime}=-\frac{1}{2}+2\epsilon$,
 $b=\frac{1}{2}+\frac{\epsilon}{2}$.  Then, the following bilinear estimate
\begin{eqnarray}
\left\|\partial_{x}(uv)\right\|_{X_{s_{1}, b^{\prime}}}\leq C\|u\|_{X_{s,b}}\|v\|_{X_{s,b}}\label{1.05}
\end{eqnarray}
fails.
\end{Theorem}
 \noindent{\bf Remark 3.} From Theorem 1.3, we know that $s_{2}>-\frac{1}{2}$ is necessary in proving Theorem 1.2.
From the proof process of Theorem 1.5,  we know that $X_{s,b}\hookrightarrow C(\R;H^{s})(b>\frac{1}{2},s\in \R)$,
which is just  Lemma 2.7 in this paper and Theorem 1.2  plays the key role  in proving Theorem 5.
 From the proof process of Theorem 1.5,  we know that $b>\frac{1}{2}$ is necessary in proving   Theorem 1.2.
Thus, we require $b>\frac{1}{2}$ in Theorem 1.3.

\begin{Theorem} \label{Theorem4}(Pointwise convergence of Kawahara equation)
Let $u_{0}\in H^{s}(\R)(s\geq\frac{1}{4})$ and $u$ be the solution to (\ref{1.01}).  Then, we have
\begin{eqnarray}
&&\lim\limits_{t\longrightarrow0}u(x,t)=u_{0}(x),\label{1.06}
\end{eqnarray}
for almost everywhere $x\in\R$.
\end{Theorem}
\noindent{\bf  Remark 4:}Follow the idea of  Proposition 4.3 of \cite{CLS}, we present the outline of Theorem 1.4. From Lemma 2.5 established in this paper, we have
\begin{eqnarray}
&&\|u\|_{L_{x}^{4}L_{t}^{\infty}}\leq C\|u\|_{X_{s,b}}(s\geq\frac{1}{4}).\label{1.07}
\end{eqnarray}
We consider the frequency truncated Kawahara equation
\begin{eqnarray}
&&\partial_{t}u_{N}+\alpha\partial_{x}^{5}u_{N}+\beta\partial_{x}^{3}u_{N}+
\partial_{x}P_{N}\left((u_{N})^{2}\right)=0,\label{1.08}\\
&& u_{N}(x,0)=P_{N}u_{0}.\label{1.09}
\end{eqnarray}
Firstly, we prove
\begin{eqnarray}
&&\lim\limits_{N\rightarrow\infty}\|u-u_{N}\|_{L_{x}^{4}L_{t}^{\infty}}=0,\label{1.010}
\end{eqnarray}
which is just Lemma 4.1 in this paper.
Since $u_{N}$ is smooth,  for all $x\in \R$,  we have
\begin{eqnarray}
&&\lim\limits_{t\rightarrow0}u_{N}(x,t)=P_{N}u_{0}(x).\label{1.011}
\end{eqnarray}
Since
\begin{eqnarray}
&&|u-u_{0}|\leq |u-u_{N}|+|u_{N}-P_{N}u_{0}|+|P_{N}u_{0}|,\label{1.012}
\end{eqnarray}
 we have
\begin{eqnarray}
&&\lim_{t\rightarrow0}\sup|u-u_{0}|\leq\lim_{t\rightarrow0}\sup|u-u_{N}|+|P^{N}u_{0}|.\label{1.013}
\end{eqnarray}
For arbitrary  $\lambda>0$,  by using the Chebyshev inequality,   (\ref{1.013})  and
 Sobolev embedding,
 we have
\begin{eqnarray}
&&|\{x\in \R:\lim\limits_{t\rightarrow0}\sup|u-u_{0}|>\lambda\}|\leq |\{x\in \R:\lim\limits_{t\rightarrow0}\sup|u-u_{N}|>\frac{\lambda}{2}\}|\nonumber\\
&&+|\{x\in \R:|P_{N}u_{0}|>\frac{\lambda}{2}\}|\leq C\lambda^{-4}\left\|u-u_{N}\right\|_{L_{x}^{4}L_{t}^{\infty}}^{4}+C\lambda^{-2}\|P^{N}u_{0}\|_{L^{2}}\nonumber\\
&&\leq C\lambda^{-4}\left\|u-u_{N}\right\|_{L_{x}^{4}L_{t}^{\infty}}^{4}+C\lambda^{-2}\|P^{N}u_{0}\|_{H^{s}}.\label{1.014}
\end{eqnarray}
Since $u_{0}\in H^{s}(\R)(s\geq\frac{1}{4})$, we have
\begin{eqnarray}
\|P^{N}u_{0}\|_{H^{s}}\rightarrow0  \label{1.015}
\end{eqnarray}
 as $N\rightarrow\infty$.  By using (\ref{1.010}) and (\ref{1.015}),  we have
\begin{eqnarray}
&&|\{x\in \R:\lim\limits_{t\rightarrow0}\sup|u-u_{0}|>\lambda\}|=0.\label{1.016}
\end{eqnarray}

\begin{Theorem} \label{Theorem5}  (Uniform convergence of Kawahara equation)
Let $u_{0}\in H^{s}(\R)(s>-\frac{1}{2})$ and  $u$ be the solution to (\ref{1.01}).  Then, we have
\begin{eqnarray}
&&\lim\limits_{t\rightarrow0}\sup\limits_{x\in\SR}|u(x,t)-U(t)u_{0}|=0.\label{1.017}
\end{eqnarray}
\end{Theorem}
\noindent{\bf  Remark 5:} Inspired by the idea of  \cite{Compaan},  we present the outline of Theorem 1.5.
Firstly, from \cite{JH2009} and the proof process of Lemma 4.1 in this paper,
 we have that
the Cauchy problem for the Kawahara equation possesses a unique solution with data in $H^{s}(\R)(s>-\frac{7}{4})$.
 From the proof process of Lemma 4.1, we have
 \begin{eqnarray}
 u-U(t)u_{0}=\eta\left(\frac{t}{T}\right)\int_{0}^{t}U(t-t^{\prime})\partial_{x}(u^{2})dt^{\prime} \label{1.018}.
 \end{eqnarray}
 From (\ref{1.018}) and Theorem 1.2, we have
\begin{eqnarray}
&&\|u-U(t)u_{0}\|_{X_{s_{1},b}}= \left\|\eta\left(\frac{t}{T}\right)\int_{0}^{t}U(t-t^{\prime})\partial_{x}(u^{2})dt^{\prime}\right\|_{X_{s_{1},b}}\nonumber\\
&&\leq C\|\partial_{x}(u^{2})\|_{X_{s_{1},b^{\prime}}}\leq C\|u\|_{X_{s_{2},b}}^{2}\leq 2C^{3}\|u_{0}\|_{H^{s_{2}}(\SR)}<\infty.\label{1.019}
\end{eqnarray}
 Here $s_{1}=\frac{1}{2}+2\epsilon, s_{2}=-\frac{1}{2}+\epsilon.$
Since $X_{s_{1},\>b}\hookrightarrow C([-T,T];H^{s_{1}}(\R))\hookrightarrow C([-T,T];C(\R))$, from (\ref{1.019}),  we have
\begin{eqnarray}
&&\lim\limits_{t\rightarrow0}\sup\limits_{x\in\SR}|u(x,t)-U(t)u_{0}|=\lim\limits_{t\rightarrow0}\sup\limits_{x\in\SR}
\left|\eta\left(\frac{t}{T}\right)\int_{0}^{t}U(t-t^{\prime})\partial_{x}(u^{2})dt^{\prime}\right|=0.\label{1.020}
\end{eqnarray}
 Here, we use Lemma 2.7 and $H^{s_{1}}(\R)\hookrightarrow  C(\R)(s_{1}=\frac{1}{2}+2\epsilon)$.

 \noindent{\bf  Remark 6:} We can use      Lemma 2.5 established in this paper and  Theorem 1.5 to
 present an alternative proof of Theorem 1.4.
  Combining Lemma 2.5 established in this paper  with the proof of Lemma 2.3 of \cite{Du},
  we immediately obtain
\begin{eqnarray}
U(t)u_{0} \longrightarrow u_{0} \qquad a.e.  \label{1.021}
\end{eqnarray}
         as $t\longrightarrow 0$ for data in $H^{s}(\R)(s\geq\frac{1}{4})$.
From (\ref{1.017}),   we know
\begin{eqnarray}
u\longrightarrow U(t)u_{0} \qquad a.e.  \label{1.022}
\end{eqnarray}
         as $t\longrightarrow 0$   for data in $H^{s}(\R)(s>-\frac{1}{2})$.
 By using the triangle inequality, we have
      \begin{eqnarray}
      \left|u-u_{0}\right|\leq \left|u-U(t)u_{0}\right|
                    + \left|u_{0}-U(t)u_{0}\right| \longrightarrow0 \qquad a.e.  \label{1.023}
      \end{eqnarray}
              as $t\longrightarrow 0$   for data in $H^{s}(\R)(s\geq\frac{1}{4})$.
 Thus, we have
  \begin{eqnarray}
  u\longrightarrow u_{0} \qquad a.e.  \label{1.022}
  \end{eqnarray}
           as $t\longrightarrow 0$   for data in $H^{s}(\R)(s\geq\frac{1}{4})$.  Thus, we give an alternative proof of Theorem 1.4.

  \noindent{\bf  Remark 7:} Compaan  et al. \cite{CLS}   studied the pointwise convergence and uniform convergence
   of the semilinear Schr\"odinger equation with rough data and random data.
   In this paper, we investigate     the pointwise convergence and uniform convergence of the Kawahara equation with rough data.
            Kawahara equation is a quasilinear evolution equation, thus,      the structure of its  is much more complicated than the structure
           of  the semilinear Schr\"odinger equation.

  \noindent{\bf  Remark 8:}  The proof of  Theorem 1.5 mainly depends on the Theorem 1.2, which is optimal  due to Theorem 1.3.
  Thus, the result of Theorem 1.5 is optimal in the sense of Theorems 1.2, 1.3.

The rest of the paper is arranged as follows. In Section 2,  we give some
preliminaries. In Section 3, we prove three bilinear estimates, which are just Theorems 1.1-1.3.
In Section 4, we give the proof of the Theorem 1.4. In Section 5, we give the proof of the Theorem 1.5.

\bigskip

\setcounter{section}{2}

\noindent{\large\bf 2. Preliminaries }

\setcounter{equation}{0}

\setcounter{Theorem}{0}

\setcounter{Lemma}{0}

\setcounter{section}{2}
In this section, we present some preliminaries.

\begin{Lemma}\label{lem2.1}
 Let  $T\in(0,1)$,  $s\in\R$,  $-\frac{1}{2}<b^{\prime}\leq 0\leq b\leq b^{\prime}+1$ and  $f\in H^{s}(\R),g\in X_{s,b^{\prime}}(\R^{2})$. Then,  we have
\begin{eqnarray}
&&\left\|\eta(t)U(t)f\right\|_{X_{s,\frac{1}{2}+\epsilon}}\leq C\left\|f\right\|_{H^{s}(\SR)},\label{2.01}\\
&&\left\|\eta\left(\frac{t}{T}\right)\int_{0}^{t}U(t-\tau)g(\tau)d\tau\right\|_{X_{s,b}(\SR^{2})}\leq CT^{1+b^{\prime}-b}\left\|g\right\|_{X_{s,b^{\prime}}(\SR^{2})}.\label{2.02}
\end{eqnarray}
\end{Lemma}

For  the proof of Lemma 2.1, we refer the readers to \cite{Bourgain1993,KPV-Duke,G2002}.
\begin{Lemma}\label{lem2.2}
Let
\begin{eqnarray*}
&&\phi(\xi)=\alpha\xi^{5}-\beta\xi^{3}, \sigma=\tau+\phi(\xi), \sigma_{j}=\tau_{j}+\phi(\xi_{j}) (1\leq j\leq 2).
\end{eqnarray*}
Then, we have
\begin{eqnarray}
&&\left|\sigma-\sigma_{1}-\sigma_{2}\right|=5|\alpha||\xi||\xi_{1}||\xi_{2}|\left|\xi^{2}+\xi_{1}^{2}-\xi\xi_{1}-\frac{3\beta}{5\alpha}\right|.\label{2.03}
\end{eqnarray}
Moreover, when $|\xi|\geq2a$ or $|\xi_{1}|\geq 2a$, where $a$ is defined as in \cite{JH2009}, then, (\ref{2.03}) implies that one of the following cases always occurs:
\begin{eqnarray}
&&\max\{|\sigma|,|\sigma_{1}|,|\sigma_{2}|\}=|\sigma|\geq C|\xi||\xi_{1}||\xi_{2}|\max\{|\xi|^{2},|\xi_{1}|^{2}\},\label{2.04}\\
&&\max\{|\sigma|,|\sigma_{1}|,|\sigma_{2}|\}=|\sigma_{1}|\geq C|\xi||\xi_{1}||\xi_{2}|\max\{|\xi|^{2},|\xi_{1}|^{2}\},\label{2.05}\\
&&\max\{|\sigma|,|\sigma_{1}|,|\sigma_{2}|\}=|\sigma_{2}|\geq C|\xi||\xi_{1}||\xi_{2}|\max\{|\xi|^{2},|\xi_{1}|^{2}\}.\label{2.06}
\end{eqnarray}
\end{Lemma}

Lemma 2.2 can be seen {\cite{JH2009}}.
\begin{Lemma}\label{lem2.3}
Let $b>\frac{1}{2}$ and $D\geq 4a$, $a$ is defined as in {\cite{JH2009}}.  Then, we have
\begin{eqnarray}
&&\|P^{D}u\|_{L_{t}^{4}L_{x}^{2}}\leq C\|u\|_{X_{0,\frac{1}{2}b}},\label{2.07}\\
&&\|P^{D}u\|_{X_{0,-\frac{b}{2}}}\leq C\|u\|_{L_{t}^{\frac{4}{3}}L_{x}^{2}},\label{2.08}\\
&&\|D_{x}^{\frac{3}{4}}P^{D}u\|_{L_{t}^{4}L_{x}^{\infty}}\leq C\|u\|_{X_{0,b}},\label{2.09}\\
&&\|u\|_{L_{xt}^{12}}\leq C\|u\|_{X_{0,b}},\label{2.010}\\
&&\|u\|_{L_{xt}^{4}}\leq C\|u\|_{X_{0,\frac{3b}{5}}}\label{2.011},\\
&&\left\|P^{D}U(t)u_{0}\right\|_{L_{x}^{4}L_{t}^{\infty}}\leq C\|u_{0}\|_{H^{\frac{1}{4}}(\SR)},\label{2.012}  \\
&&\left\|P^{D}D_{x}^{\frac{3}{8}}u\right\|_{L_{xt}^{4}}\leq C\|u\|_{X_{0,b}}. \label{2.013}
\end{eqnarray}
\end{Lemma}
\noindent {\bf Proof.}
For the proof of  (\ref{2.07})-(\ref{2.09}), we refer the readers to  Lemma 2.6 of  {\cite{YL}}. For the proof of (\ref{2.010}), we refer the readers to
(2.21) of {\cite{JH2009}}. Interpolating (\ref{2.010}) with
\begin{eqnarray*}
\|u\|_{L_{xt}^{2}}=C\|u\|_{L_{\xi\tau}^{2}}
\end{eqnarray*}
yields (\ref{2.011}).    For the proof of (\ref{2.012}), (\ref{2.013}),  we refer the readers to (2.10) and (2.13) of \cite{JH2009}, respectively.

We have completed the proof of Lemma 2.3.

\begin{Lemma}\label{lem2.4}
Let  $s\geq\frac{1}{4}$.  Then we have
\begin{eqnarray}
&&\|U(t)u_{0}\|_{L_{x}^{4}L_{t}^{\infty}}\leq C\|u_{0}\|_{H^{s}(\SR)}\label{2.014}.
\end{eqnarray}
\end{Lemma}
\noindent {\bf Proof.}By using (\ref{2.012}) and the Sobolev embeddings Theorem $W^{\frac{1}{4}+\epsilon,4}(\R)$, we have
\begin{eqnarray}
&&\|U(t)u_{0}\|_{L_{x}^{4}L_{t}^{\infty}}\leq C\left\|P^{D}U(t)u_{0}\right\|_{L_{x}^{4}L_{t}^{\infty}}+C\left\|P_{D}U(t)u_{0}\right\|_{L_{x}^{4}L_{t}^{\infty}}\nonumber\\
&&\leq C\|u_{0}\|_{H^{\frac{1}{4}}(\SR)}+C\left\|D_{t}^{\frac{1}{4}+\epsilon}P_{D}U(t)u_{0}\right\|_{L_{xt}^{4}}\nonumber\\
&&\leq C\|u_{0}\|_{H^{\frac{1}{4}}(\SR)}+
C\left\|U(t)\mathscr{F}_{x}^{-1}\left(\left|\alpha\xi^{5}+\beta\xi^{3}\right|^{\frac{1}{4}+\epsilon}\chi_{|\xi|\leq D}(\xi)\mathscr{F}_{x}u_{0}(\xi)\right)\right\|_{L_{xt}^{4}}\nonumber\\
&&\leq C\|u_{0}\|_{H^{\frac{1}{4}}(\SR)}+
C\left\|\mathscr{F}_{x}^{-1}\left(\left|\alpha\xi^{5}+\beta\xi^{3}\right|^{\frac{1}{4}+\epsilon}\chi_{|\xi|\leq D}(\xi)\mathscr{F}_{x}u_{0}(\xi)\right)\right\|_{L_{xt}^{2}}\nonumber\\
&&\leq C\|u_{0}\|_{H^{\frac{1}{4}}(\SR)}+
C\left\|\left|\alpha\xi^{5}+\beta\xi^{3}\right|^{\frac{1}{4}+\epsilon}\chi_{|\xi|\leq D}(\xi)\mathscr{F}_{x}u_{0}(\xi)\right\|_{L_{\xi\tau}^{2}}\nonumber\\
&&\leq C\|u_{0}\|_{H^{\frac{1}{4}}(\SR)}\leq C\|u_{0}\|_{H^{s}(\SR)}.\label{2.015}
\end{eqnarray}
 Here,     $a$ is defined as in {\cite{JH2009}} and $D\geq4a.$

We have completed the proof of Lemma 2.4.

\begin{Lemma}\label{lem2.5}
Let $s\geq\frac{1}{4}$ and $b>\frac{1}{2}$. Then, we have
\begin{eqnarray}
&&\|u\|_{L_{x}^{4}L_{t}^{\infty}}\leq C\|u\|_{X_{s,b}}.\label{2.016}
\end{eqnarray}
\end{Lemma}
\noindent {\bf Proof.}
By changing variable $\tau=\lambda-\phi(\xi)$, we derive
 \begin{eqnarray}
&&u(x,t)=\frac{1}{2\pi}\int_{\SR^{2}}e^{ix\xi+it\tau}\mathscr{F}u(\xi,\tau)d\xi d\tau\nonumber\\
&&=\frac{1}{2\pi}\int_{\SR^{2}}e^{ix\xi+it(\lambda-\phi(\xi))}\mathscr{F}u(\xi,\lambda-\phi(\xi))d\xi d\lambda\nonumber\\
&&=\frac{1}{2\pi}\int_{\SR}e^{it\lambda}\left(\int_{\SR}e^{ix\xi-it\phi(\xi)}
\mathscr{F}u(\xi,\lambda-\phi(\xi)) d\xi\right)d\lambda\label{2.017}.
 \end{eqnarray}
 By using  (\ref{2.014}),  (\ref{2.017}) and Minkowski's inequality,  for  $b>\frac{1}{2},$  we derive
\begin{eqnarray}
&&\left\|u\right\|_{L_{x}^{4}L_{t}^{\infty}}\leq C\int_{\SR}\left\|\left(\int_{\SR}e^{ix\xi-it\phi(\xi)}
\mathscr{F}u(\xi,\lambda-\phi(\xi))d\xi\right)\right\|_{L_{x}^{4}L_{t}^{\infty}} d\lambda\nonumber\\&&\leq C\int_{\SR}\left\|\mathscr{F}u(\xi,\lambda-\phi(\xi))\right\|_{H^{s}}d\lambda\nonumber\\&&\leq C
\left[\int_{\SR}(1+|\lambda|)^{2b}\left\|\mathscr{F}u(\xi,\lambda-\phi(\xi))\right\|_{H^{s}}^{2}d\lambda\right]^{\frac{1}{2}}
\left[\int_{\SR}(1+|\lambda|)^{-2b}d\lambda\right]^{\frac{1}{2}}\nonumber\\
&&\leq C\left[\int_{\SR}(1+|\tau+\phi(\xi)|)^{2b}\left\|\mathscr{F}u(\xi,\tau)\right\|_{H^{s}}^{2}d\tau\right]^{\frac{1}{2}}
=\left\|u\right\|_{X_{s,b}}.\label{2.018}
\end{eqnarray}

This completes the proof of Lemma 2.5.

\begin{Lemma}\label{lem2.6}
Let $b=\frac{1}{2}+\epsilon$. Then, we have
\begin{eqnarray}
&&\|I(u_{1},u_{2})\|_{L_{xt}^{2}}\leq C\prod_{j=1}^{2}\|u_{j}\|_{X_{0,b}},\label{2.018}\\
&&\|I_{k}(u_{1},u_{2})\|_{L_{xt}^{2}}\leq C\prod_{j=1}^{2}\|u_{j}\|_{X_{0,b}} \qquad (1\leq k\leq 2),\label{2.019}
\end{eqnarray}
where
\begin{eqnarray*}
&&\mathscr{F}I(u_{1},u_{2})(\xi,\tau)=\int_{\xi=\xi_{1}+\xi_{2},|\xi|\geq 4a,
\tau=\tau_{1}+\tau_{2}}|\xi_{1}^{4}-\xi_{2}^{4}|^{\frac{1}{2}}\mathscr{F} u_{1}(\xi_{1},\tau_{1})\mathscr{F}u_{2}(\xi_{2},\tau_{2})d\xi_{1}d\tau_{1},\\
&&\mathscr{F}I_{k}(u_{1},u_{2})(\xi,\tau)=\int_{\xi=\xi_{1}+\xi_{2},|\xi_{k}|\geq 4a,
\tau=\tau_{1}+\tau_{2}}|\xi_{1}^{4}-\xi_{2}^{4}|^{\frac{1}{2}}\prod\limits_{j=1}^{2}\mathscr{F} u_{j}(\xi_{j},\tau_{j})d\xi_{1}d\tau_{1}  \quad (1\leq k\leq2).
\end{eqnarray*}
\end{Lemma}

For the proof of Lemma 2.6, we refer the readers to  Theorem  3.1 of  \cite{YL}.

\begin{Lemma}\label{lem2.7}
Let $b>\frac{1}{2}$. Then, we have $X_{s,b}(\R^{2})\hookrightarrow C(\R;H^{s}(\R))$.
\end{Lemma}

For the proof of Lemma 2.7, we refer the readers to  Lemma 4 of  \cite{ET2013}.

\bigskip

\setcounter{section}{3}

\noindent{\large\bf 3. Bilinear estimates }

\setcounter{equation}{0}

\setcounter{Theorem}{0}

\setcounter{Lemma}{0}

\setcounter{section}{3}
In this section, we prove Theorems 1.1-1.3.

To prove Theorem 1.1, it suffices to prove Lemma 3.1.

\begin{Lemma}\label{lem3.1}
Let  $s\geq-\frac{7}{4}+4\epsilon$,  $b^{\prime}=-\frac{1}{2}+2\epsilon$ and $b=\frac{1}{2}+\frac{\epsilon}{2}$. Then, we have
\begin{eqnarray}
&&\left\|\partial_{x}\left(u_{1}u_{2}\right)\right\|_{X_{s,b^{\prime}}}\leq C\prod\limits_{j=1}^{2}\|u_{j}\|_{X_{s,b}}.\label{3.01}
\end{eqnarray}
\end{Lemma}
\noindent {\bf Proof.}
To prove (\ref{3.01}), by duality, it suffices to prove
\begin{eqnarray}
&&\left|\int_{\SR^{2}}\partial_{x}(u_{1}u_{2})\bar{h}dxdt\right|\leq C\|h\|_{X_{-s,-b^{\prime}}}\prod\limits_{j=1}^{2}\|u_{j}\|_{X_{s,b}}.\label{3.02}
\end{eqnarray}
We define
\begin{eqnarray*}
&&\int_{*}=\int_{\xi=\xi_{1}+\xi_{2},\tau=\tau_{1}+\tau_{2}},\\
&&f_{j}(\xi_{j},\tau_{j})=\langle\xi_{j}\rangle^{s}\langle\sigma_{j}\rangle^{b}\mathscr{F}u_{j}(\xi_{j},\tau_{j})(j=1,2),\\
&&g(\xi,\tau)=\langle\xi\rangle^{-s}\langle\sigma\rangle^{-b^{\prime}}\mathscr{F}h(\xi,\tau).
\end{eqnarray*}
To prove (\ref{3.02}),    it suffices to prove
\begin{eqnarray}
&&\int_{\SR^{2}}\int_{*}\frac{|\xi|\langle\xi\rangle^{s}}
{\langle\sigma\rangle^{-b^{\prime}}\prod\limits_{j=1}^{2}\langle\sigma_{j}\rangle^{b}
\prod\limits_{j=1}^{2}\langle\xi_{j}\rangle^{s}}f_{1}(\xi_{1},\tau_{1})f_{2}(\xi_{2},\tau_{2})g(\xi,\tau)d\xi_{1}d\tau_{1}d\xi d\tau\nonumber\\
&&\leq C\|g\|_{L_{\xi\tau}^{2}}\|f_{1}\|_{L_{\xi\tau}^{2}}\|f_{2}\|_{L_{\xi\tau}^{2}} \label{3.03}
\end{eqnarray}
with  the aid of   the Plancherel identity.
We define
\begin{eqnarray*}
&&K_{1}(\xi_{1},\tau_{1},\xi,\tau)=\frac{|\xi|\langle\xi\rangle^{s}}
{\langle\sigma\rangle^{-b^{\prime}}\prod\limits_{j=1}^{2}\langle\sigma_{j}\rangle^{b}\prod\limits_{j=1}^{2}\langle\xi_{j}\rangle^{s}},\\
&&\mathscr{F}F_{j}(\xi_{j},\tau_{j})=\frac{f_{j}(\xi_{j},\tau_{j})}{\langle\sigma_{j}\rangle^{b}}(j=1,2),\mathscr{F}G(\xi,\tau)=\frac{g(\xi,\tau)}{\langle\sigma\rangle^{-b^{\prime}}},\\
&&I_{1}=\int_{\SR^{2}}\int_{*}K_{1}(\xi_{1},\tau_{1},\xi,\tau)f_{1}(\xi_{1},\tau_{1})f_{2}(\xi_{2},\tau_{2})g(\xi,\tau)d\xi_{1}d\tau_{1}d\xi d\tau.
\end{eqnarray*}
Without loss of generality, we can assume that $|\xi_{1}|\geq |\xi_{2}|$. Obviously,
\begin{eqnarray*}
&&\Omega=\left\{(\xi_{1},\xi,\tau_{1},\tau)\in\R^{4}:\xi=\xi_{1}+\xi_{2},\tau=\tau_{1}+\tau_{2},|\xi_{1}|\geq|\xi_{2}|\right\}\subset\bigcup_{j=1}^{6}\Omega_{j},
\end{eqnarray*}
and
\begin{eqnarray*}
&&\Omega_{1}=\left\{(\xi_{1},\tau_{1},\xi,\tau)\in \Omega:|\xi_{1}|\leq 4a\right\},\\
&&\Omega_{2}=\left\{(\xi_{1},\tau_{1},\xi,\tau)\in \Omega:|\xi_{1}|> 4a, |\xi_{1}|>4|\xi_{2}|,|\xi_{2}|\leq a\right\},\\
&&\Omega_{3}=\left\{(\xi_{1},\tau_{1},\xi,\tau)\in \Omega:|\xi_{1}|> 4a, |\xi_{1}|>4|\xi_{2}|,|\xi_{2}|> a\right\},\\
&&\Omega_{4}=\left\{(\xi_{1},\tau_{1},\xi,\tau)\in \Omega:|\xi_{1}|> 4a,|\xi_{2}|\leq |\xi_{1}|\leq 4|\xi_{2}|,\xi_{1}\xi_{2}\geq0\right\},\\
&&\Omega_{5}=\left\{(\xi_{1},\tau_{1},\xi,\tau)\in \Omega:|\xi_{1}|> 4a,|\xi_{2}|\leq |\xi_{1}|\leq 4|\xi_{2}|,\xi_{1}\xi_{2}<0,4|\xi|\geq|\xi_{2}|\right\},\\
&&\Omega_{6}=\left\{(\xi_{1},\tau_{1},\xi,\tau)\in \Omega:|\xi_{1}|> 4a,|\xi_{2}|\leq |\xi_{1}|\leq 4|\xi_{2}|,\xi_{1}\xi_{2}<0,4|\xi|<|\xi_{2}|\right\}.
\end{eqnarray*}
(1) When $(\xi_{1},\xi,\tau_{1},\tau)\in \Omega_{1}$, which yield $|\xi|\leq |\xi_{1}|+|\xi_{2}|\leq 8a$, therefore, we have
\begin{eqnarray}
&&K_{1}(\xi_{1},\tau_{1},\xi,\tau)\leq \frac{C}{\langle\sigma\rangle^{-b^{\prime}}
\prod\limits_{j=1}^{2}\langle\sigma_{j}\rangle^{b}}\leq \frac{C}{\prod\limits_{j=1}^{2}\langle\sigma_{j}\rangle^{b}}.\label{3.04}
\end{eqnarray}
By using (\ref{3.04}),  the Cauchy-Schwarz inequality,  the Plancherel identity and the H\"older inequality as well as    (\ref{2.011}),  we have
\begin{eqnarray}
&&I_{1}\leq C\int_{\SR^{2}}\int_{*}\frac{f_{1}f_{2}g}{\prod\limits_{j=1}^{2}\langle\sigma_{j}\rangle^{b}}d\xi_{1}d\tau_{1}d\xi d\tau\nonumber\\
&&\leq C\left\|\int_{*}\frac{f_{1}f_{2}}{\prod\limits_{j=1}^{2}\langle\sigma_{j}\rangle^{b}}d\xi_{1}d\tau_{1}\right\|_{L_{\xi\tau}^{2}}\|g\|_{L_{\xi\tau}^{2}}\nonumber\\
&&\leq C\|F_{1}F_{2}\|_{L_{xt}^{2}}\|g\|_{L_{\xi\tau}^{2}}\nonumber\\
&&\leq C\|F_{1}\|_{L_{xt}^{4}}\|F_{2}\|_{L_{xt}^{4}}\|g\|_{L_{\xi\tau}^{2}}\nonumber\\&&
\leq C\|F_{1}\|_{X_{0,\frac{3b}{5}}}\|F_{2}\|_{X_{0,\frac{3b}{5}}}\nonumber\\&&
\leq C\|f_{1}\|_{L_{\xi\tau}^{2}}\|f_{2}\|_{L_{\xi\tau}^{2}}\|g\|_{L_{\xi\tau}^{2}}.\label{3.05}
\end{eqnarray}
(2)When $(\xi_{1},\xi,\tau_{1},\tau)\in \Omega_{2}$, which yield $|\xi_{1}|\sim|\xi|$ and $|\xi_{2}|\leq a$, therefore, we have
\begin{eqnarray}
&&K_{1}(\xi_{1},\tau_{1},\xi,\tau)\leq C\frac{|\xi|}
{\langle\sigma\rangle^{-b^{\prime}}\prod\limits_{j=1}^{2}\langle\sigma_{j}\rangle^{b}}\leq C
\frac{|\xi_{1}^{4}-\xi_{2}^{4}|^{\frac{1}{2}}}{\prod\limits_{j=1}^{2}\langle\sigma_{j}\rangle^{b}}.\label{3.06}
\end{eqnarray}
By using (\ref{3.06}), the Cauchy-Schwarz inequality and the Plancherel identity as well as Lemma 2.6, we have
\begin{eqnarray}
&&I_{1}\leq C\int_{\SR^{2}}\int_{*}\frac{|\xi_{1}^{4}-\xi_{2}^{4}|^{\frac{1}{2}}}
{\prod\limits_{j=1}^{2}\langle\sigma_{j}\rangle^{b}}f_{1}f_{2}gd\xi_{1}d\tau_{1}d\xi d\tau\nonumber\\
&&\leq C\left\|\int_{*}\frac{|\xi_{1}^{4}-\xi_{2}^{4}|^{\frac{1}{2}}}{\prod\limits_{j=1}^{2}
\langle\sigma_{j}\rangle^{b}}f_{1}f_{2}d\xi_{1}d\tau_{1}\right\|_{L_{\xi\tau}^{2}}\|g\|_{L_{\xi\tau}^{2}}\nonumber\\
&&\leq C\|F_{1}\|_{X_{0,b}}\|F_{2}\|_{X_{0,b}}\|g\|_{L_{\xi\tau}^{2}}\leq C
\|f_{1}\|_{L_{\xi\tau}^{2}}\|f_{2}\|_{L_{\xi\tau}^{2}}\|g\|_{L_{\xi\tau}^{2}}.\label{3.07}
\end{eqnarray}
(3)When$(\xi_{1},\xi,\tau_{1},\tau)\in \Omega_{3}$, which yield $|\xi_{1}|\sim|\xi|,|\xi_{1}|\geq4a>2a$, $|\xi_{2}|>a$,
 then, we consider (\ref{2.04})-(\ref{2.06}),  respectively.\\
When (\ref{2.04}) is valid, since $s\geq-\frac{7}{4}+4\epsilon$, then, we have
\begin{eqnarray}
&&K_{1}(\xi_{1},\tau_{1},\xi,\tau)\leq \frac{|\xi_{1}|^{1+4b^{\prime}}|\xi_{2}|^{b^{\prime}-s}}
{\prod\limits_{j=1}^{2}\langle\sigma_{j}\rangle^{b}}
\leq C\frac{|\xi_{1}|^{\frac{3}{4}+6\epsilon}}{\prod\limits_{j=1}^{2}\langle\sigma_{j}\rangle^{b}}\leq C
\frac{|\xi_{1}^{4}-\xi_{2}^{4}|^{\frac{1}{2}}}{\prod\limits_{j=1}^{2}\langle\sigma_{j}\rangle^{b}}.\label{3.08}
\end{eqnarray}
This case can be proved similarly to Case (2).\\
When (\ref{2.05}) is valid, which yield $\langle\sigma_{1}\rangle^{-b}\langle\sigma\rangle^{b^{\prime}}\leq
\langle\sigma_{1}\rangle^{b^{\prime}}\langle\sigma\rangle^{-b}$, since $s\geq-\frac{7}{4}+4\epsilon$,  we have
\begin{eqnarray}
&&K_{1}(\xi_{1},\tau_{1},\xi,\tau)\leq \frac{|\xi|^{1+4b^{\prime}}|\xi_{2}|^{-s+b^{\prime}}}
{\langle\sigma_{2}\rangle^{b}\langle\sigma\rangle^{b}}
\leq C\frac{|\xi_{1}|^{\frac{1}{4}+6\epsilon}}{\langle\sigma_{2}\rangle^{b}\langle\sigma\rangle^{b}}\leq C
\frac{|\xi^{4}-\xi_{2}^{4}|^{\frac{1}{2}}}{\langle\sigma_{2}\rangle^{b}\langle\sigma\rangle^{b}}.\label{3.09}
\end{eqnarray}
This case can be proved similarly to Case (2).\\
When (\ref{2.06}) is valid, which yield $\langle\sigma_{2}\rangle^{-b}\langle\sigma\rangle^{b^{\prime}}\leq
 \langle\sigma_{2}\rangle^{b^{\prime}}\langle\sigma\rangle^{-b}$, since $s\geq-\frac{7}{4}+4\epsilon$, then we have
\begin{eqnarray}
&&K_{1}(\xi_{1},\tau_{1},\xi,\tau)\leq \frac{|\xi|^{1+4b^{\prime}}|\xi_{2}|^{-s+b^{\prime}}}
{\langle\sigma_{1}\rangle^{b}\langle\sigma\rangle^{b}}
\leq C\frac{|\xi_{1}|^{\frac{1}{4}+6\epsilon}}{\langle\sigma_{1}\rangle^{b}\langle\sigma\rangle^{b}}\leq C
\frac{|\xi_{1}|^{\frac{3}{4}}}{\langle\sigma_{1}\rangle^{b}\langle\sigma\rangle^{b}}.\label{3.010}
\end{eqnarray}
By using (\ref{3.010}),  the Cauchy-Schwarz inequality,  (\ref{2.08})-(\ref{2.09}) and   the H\"older inequality, we have
\begin{eqnarray}
&&I\leq C\int_{\SR^{2}}\int_{*}\frac{|\xi_{1}|^{\frac{3}{4}}}{\langle\sigma_{1}\rangle^{b}\langle\sigma\rangle^{b}}f_{1}f_{2}gd\xi_{1}d\tau_{1}d\xi d\tau\nonumber\\
&&\leq C\left\|\langle\sigma\rangle^{-b}\int_{*}\frac{|\xi_{1}|^{\frac{3}{4}}}{\langle\sigma_{1}\rangle^{b}}f_{1}f_{2}d\xi_{1}d\tau_{1}\right\|
_{L_{\xi\tau}^{2}}\|h\|_{L_{\xi\tau}^{2}}\leq C\left\|P^{4a}((D_{x}^{\frac{3}{4}}P^{D}F_{1})\mathscr{F}^{-1}f_{2})\right\|_{X_{0,-b}}\|g\|_{L_{\xi\tau}^{2}}\nonumber\\
&&\leq C\|D_{x}^{\frac{3}{4}}P^{D}F_{1}\mathscr{F}^{-1}f_{2}\|_{L_{t}^{\frac{4}{3}}L_{x}^{2}}\|g\|_{L_{\xi\tau}^{2}}\leq C\left\|D_{x}^{\frac{3}{4}}P^{D}F_{1}\right\|_{L_{t}^{4}L_{x}^{\infty}}\|\mathscr{F}^{-1}f_{2}\|_{L_{xt}^{2}}\|g\|_{L_{\xi\tau}^{2}}\nonumber\\
&&\leq C\|f_{1}\|_{L_{\xi\tau}^{2}}\|f_{2}\|_{L_{\xi\tau}^{2}}\|g\|_{L_{\xi\tau}^{2}}.\label{3.011}
\end{eqnarray}
(4)When  $(\xi_{1},\xi,\tau_{1},\tau)\in \Omega_{4}$, which yield $|\xi_{1}|\sim|\xi_{2}|\sim|\xi|,|\xi_{1}|\geq4a$,
 then, we consider (\ref{2.04})-(\ref{2.06}), respectively.\\
When (\ref{2.04}) is valid,   since  $s\geq-\frac{7}{4}+4\epsilon$, we have
\begin{eqnarray}
&&K_{1}(\xi_{1},\tau_{1},\xi,\tau)\leq C\frac{|\xi_{1}|^{1+5b^{\prime}-s}}
{\prod\limits_{j=1}^{2}\langle\sigma_{j}\rangle^{b}}
\leq C\frac{|\xi_{1}|^{\frac{1}{4}+6\epsilon}}{\prod\limits_{j=1}^{2}\langle\sigma_{j}\rangle^{b}}\leq C
\frac{|\xi_{1}|^{\frac{3}{4}}}{\prod\limits_{j=1}^{2}\langle\sigma_{j}\rangle^{b}}.\label{3.012}
\end{eqnarray}
By using (\ref{3.012}), the Cauchy-Schwarz inequality,
 and the Plancherel identity as well as (\ref{2.07}),  (\ref{2.08}), we have
\begin{eqnarray}
&&I_{1}\leq C\int_{\SR^{2}}\int_{*}\frac{|\xi_{1}|^{\frac{3}{4}}}{\prod\limits_{j=1}^{2}\langle\sigma_{j}\rangle^{b}}f_{1}f_{2}gd\xi_{1}d\tau_{1}d\xi d\tau\nonumber\\
&&\leq C\left\|\int_{\SR^{2}}\int_{*}\frac{\prod\limits_{j=1}^{2}|\xi_{j}|^{\frac{3}{8}}}{\prod\limits_{j=1}^{2}\langle\sigma_{j}\rangle^{b}}f_{1}f_{2}d\xi_{1}d\tau_{1}\right\|_{L_{\xi\tau}^{2}}\|g\|_{L_{\xi\tau}^{2}}\leq C\|g\|_{L_{\xi\tau}^{2}}\prod\limits_{j=1}^{2}\|D_{x}^{\frac{3}{8}}P^{D}F_{j}\|_{L_{tx}^{4}}\nonumber\\
&&\leq C\|F_{1}\|_{X_{0,b}}\|F_{2}\|_{X_{0,\frac{b}{2}}}\|g\|_{L_{\xi\tau}^{2}}\leq C\|f_{1}\|_{L_{\xi\tau}^{2}}\|f_{2}\|_{L_{\xi\tau}^{2}}\|g\|_{L_{\xi\tau}^{2}}.\label{3.013}
\end{eqnarray}
When (\ref{2.05}) is valid, which yield $\langle\sigma_{1}\rangle^{-b}\langle\sigma\rangle^{b^{\prime}}\leq
 \langle\sigma_{1}\rangle^{b^{\prime}}\langle\sigma\rangle^{-b}$,  since $s\geq-\frac{7}{4}+4\epsilon$, we have
\begin{eqnarray}
&&K_{1}(\xi_{1},\tau_{1},\xi,\tau)\leq C\frac{|\xi_{1}|^{1+5b^{\prime}-s}}
{\langle\sigma_{2}\rangle^{b}\langle\sigma\rangle^{b}}
\leq C\frac{|\xi_{2}|^{\frac{1}{4}+6\epsilon}}{\langle\sigma_{2}\rangle^{b}\langle\sigma\rangle^{b}}\leq
C\frac{|\xi_{2}|^{\frac{3}{8}}|\xi|^{\frac{3}{8}}}{\langle\sigma_{2}\rangle^{b}\langle\sigma\rangle^{b}}.\label{3.014}
\end{eqnarray}
This case can be proved similarly to (\ref{3.013}).\\
When (\ref{2.06}) is valid, which yield $\langle\sigma_{2}\rangle^{-b}\langle\sigma\rangle^{b^{\prime}}\leq
 \langle\sigma_{2}\rangle^{b^{\prime}}\langle\sigma\rangle^{-b}$,  since $s\geq-\frac{7}{4}+4\epsilon$, we have
\begin{eqnarray}
&&K_{1}(\xi_{1},\tau_{1},\xi,\tau)\leq C\frac{|\xi_{1}|^{1+5b^{\prime}-s}}
{\langle\sigma_{1}\rangle^{b}\langle\sigma\rangle^{b}}
\leq C\frac{|\xi_{1}|^{\frac{1}{4}+6\epsilon}}{\langle\sigma_{1}\rangle^{b}\langle\sigma\rangle^{b}}\leq
C\frac{|\xi_{1}|^{\frac{3}{8}}|\xi|^{\frac{3}{8}}}{\langle\sigma_{1}\rangle^{b}\langle\sigma\rangle^{b}}.\label{3.015}
\end{eqnarray}
This case can be proved similarly to (\ref{3.013}).\\

\noindent(5)When$(\xi_{1},\xi,\tau_{1},\tau)\in \Omega_{5}$, which yield $|\xi_{1}|\sim|\xi_{2}|\sim|\xi|,|\xi_{1}|\geq 4a$.\\

This case can be proved similarly to Case (4).\\

\noindent(6)When$(\xi_{1},\xi,\tau_{1},\tau)\in \Omega_{6}$, which yield $|\xi_{1}|\sim|\xi_{2}|,|\xi_{1}|\geq4a,4|\xi|<|\xi_{2}|$.
Then, we consider (\ref{2.04})-(\ref{2.06}),   respectively.\\

\noindent When (\ref{2.04}) is valid, since $s\geq-\frac{7}{4}+4\epsilon$,  we have
\begin{eqnarray}
&&K_{1}(\xi_{1},\tau_{1},\xi,\tau)\leq C\frac{|\xi|^{1+b^{\prime}}|\xi_{1}|^{4b^{\prime}-2s}}
{\prod\limits_{j=1}^{2}\langle\sigma_{j}\rangle^{b}}
\leq C\frac{|\xi|^{\frac{1}{2}}|\xi_{1}|^{\frac{3}{2}}}{\prod\limits_{j=1}^{2}\langle\sigma_{j}\rangle^{b}}\leq C
 \frac{|\xi_{1}^{2}-\xi_{2}^{2}|^{\frac{1}{2}}}{\prod\limits_{j=1}^{2}\langle\sigma_{j}\rangle^{b}}    .\label{3.016}
\end{eqnarray}
This case can be proved similarly to (\ref{3.07}).\\
When (\ref{2.05}) is valid, which yield $\langle\sigma_{1}\rangle^{-b}\langle\sigma\rangle^{b^{\prime}}\leq
\langle\sigma_{1}\rangle^{-b^{\prime}}\langle\sigma\rangle^{-b}$, since $s\geq-\frac{7}{4}+4\epsilon$,  we have
\begin{eqnarray}
&&K_{1}(\xi_{1},\tau_{1},\xi,\tau)\leq C\frac{|\xi|^{1+b^{\prime}}|\xi_{1}|^{4b^{\prime}-2s}}
{\langle\sigma\rangle^{b}\langle\sigma_{2}\rangle^{b}}
\leq C\frac{|\xi|^{\frac{1}{2}}|\xi_{1}|^{\frac{3}{2}}}{\langle\sigma\rangle^{b}\langle\sigma_{2}\rangle^{b}}\leq C
\frac{|\xi^{2}-\xi_{2}^{2}|^{\frac{1}{2}}}{\langle\sigma\rangle^{b}\langle\sigma_{2}\rangle^{b}}.\label{3.017}
\end{eqnarray}
This case can be proved similarly to (\ref{3.07}).\\
When (\ref{2.06}) is valid, which yield $\langle\sigma_{2}\rangle^{-b}\langle\sigma\rangle^{b^{\prime}}\leq
\langle\sigma_{2}\rangle^{b^{\prime}}\langle\sigma\rangle^{-b}$, since $s\geq-\frac{7}{4}+4\epsilon$, we have
\begin{eqnarray}
&&K_{1}(\xi_{1},\tau_{1},\xi,\tau)\leq C\frac{|\xi|^{1+b^{\prime}}|\xi_{1}|^{4b^{\prime}-2s}}
{\langle\sigma\rangle^{b}\langle\sigma_{1}\rangle^{b}}
\leq C\frac{|\xi|^{\frac{1}{2}}|\xi_{1}|^{\frac{3}{2}}}{\langle\sigma\rangle^{b}\langle\sigma_{1}\rangle^{b}}\leq C
\frac{|\xi^{2}-\xi_{1}^{2}|^{\frac{1}{2}}}{\langle\sigma\rangle^{b}\langle\sigma_{1}\rangle^{b}}.\label{3.018}
\end{eqnarray}
This case can be proved similarly to (\ref{3.07}).

This completes the proof of Lemma 3.1.

To prove Theorem 1.2, it suffices to prove Lemma 3.2.

\begin{Lemma}\label{lem3.2}
Let $s_{1}=\frac{1}{2}+2\epsilon$, $b=\frac{1}{2}+\frac{\epsilon}{2}$, $b^{\prime}=-\frac{1}{2}+2\epsilon,$  $s_{2}\geq-\frac{1}{2}+\epsilon$.  Then,  we have
\begin{eqnarray}
&&\|\partial_{x}(u_{1}u_{2})\|_{X_{s_{1},b^{\prime}}}\leq C\|u_{1}\|_{X_{s_{2},b}}\|u_{2}\|_{X_{s_{2},b}}.\label{3.019}
\end{eqnarray}
\end{Lemma}
\noindent {\bf Proof.}
To prove (\ref{3.01}), by duality, it suffices to prove
\begin{eqnarray}
&&\left|\int_{\SR^{2}}\partial_{x}(u_{1}u_{2})\bar{h}dxdt\right|\leq C\|h\|_{X_{-s_{1},-b^{\prime}}}\|u_{1}\|_{X_{s_{2},b}}\|u_{2}\|_{X_{s_{2},b}}.\label{3.020}
\end{eqnarray}
We define
\begin{eqnarray*}
&&\int_{*}=\int_{\xi=\xi_{1}+\xi_{2},\tau=\tau_{1}+\tau_{2}},\\
&&f_{j}(\xi_{j},\tau_{j})=\langle\xi_{j}\rangle^{s_{2}}\langle\sigma_{j}\rangle^{b}\mathscr{F}u_{j}(\xi_{j},\tau_{j})(j=1,2),\\
&&g(\xi,\tau)=\langle\xi\rangle^{-s_{1}}\langle\sigma\rangle^{-b^{\prime}}\mathscr{F}h(\xi,\tau).
\end{eqnarray*}
 To prove (\ref{3.020}),     it suffices to prove
\begin{eqnarray}
&&\int_{\SR^{2}}\int_{*}\frac{|\xi|\langle\xi\rangle^{s_{1}}}
{\langle\sigma\rangle^{-b^{\prime}}\prod\limits_{j=1}^{2}\langle\xi_{j}\rangle^{s_{2}}\langle\sigma_{j}\rangle^{b}}f_{1}(\xi_{1},\tau_{1})f_{2}(\xi_{2},\tau_{2})g(\xi,\tau)d\xi_{1}d\tau_{1}d\xi d\tau\nonumber\\
&&\leq C\|g\|_{L_{\xi\tau}^{2}}\|f_{1}\|_{L_{\xi\tau}^{2}}\|f_{2}\|_{L_{\xi\tau}^{2}}.\label{3.021}
\end{eqnarray}
with the aid of   the Plancherel identity.
We define
\begin{eqnarray*}
&&K_{2}(\xi_{1},\tau_{1},\xi,\tau)=\frac{|\xi|\langle\xi\rangle^{s_{1}}}
{\langle\sigma\rangle^{-b^{\prime}}\prod\limits_{j=1}^{2}\langle\xi_{j}\rangle^{s_{2}}\langle\sigma_{j}\rangle^{b}},\\
&&\mathscr{F}F_{j}(\xi_{j},\tau_{j})=\frac{f_{j}(\xi_{j},\tau_{j})}{\langle\sigma_{j}\rangle^{b}}(j=1,2),\mathscr{F}G(\xi,\tau)=\frac{g(\xi,\tau)}{\langle\sigma\rangle^{-b^{\prime}}},\\
&&I=\int_{\SR^{2}}\int_{*}K_{2}(\xi_{1},\tau_{1},\xi_{2},\tau_{2})f_{1}(\xi_{1},\tau_{1})f_{2}(\xi_{2},\tau_{2})g(\xi,\tau)d\xi_{1}d\tau_{1}d\xi d\tau.
\end{eqnarray*}

Without loss of generality, we can assume that $|\xi_{1}|\geq |\xi_{2}|$, we have
\begin{eqnarray*}
&&A=\left\{(\xi_{1},\xi,\tau_{1},\tau)\in\R^{4}:\xi=\xi_{1}+\xi_{2},\tau=\tau_{1}+\tau_{2},|\xi_{1}|\geq|\xi_{2}|\right\}\subset\bigcup_{j=1}^{6}A_{j},
\end{eqnarray*}
where $A_{j}(1\leq j\leq6,j\in N)$ are defined as in Lemma 3.1.

\noindent (1) When $(\xi_{1},\xi,\tau_{1},\tau)\in A_{1}$, which yield $|\xi|\leq |\xi_{1}|+|\xi_{2}|\leq 8a$, therefore, we have
\begin{eqnarray}
&&K_{2}(\xi_{1},\tau_{1},\xi,\tau)\leq \frac{C}{\langle\sigma\rangle^{-b^{\prime}}\prod\limits_{j=1}^{2}\langle\sigma_{j}\rangle^{b}}.\label{3.022}
\end{eqnarray}

This case can be proved similarly to Case (1) of Lemma 3.1.

\noindent (2)When $(\xi_{1},\xi,\tau_{1},\tau)\in A_{2}$, which yield $|\xi_{1}|\sim|\xi|$ and $|\xi_{2}|\leq a$, therefore, we have
\begin{eqnarray}
&&K_{2}(\xi_{1},\tau_{1},\xi,\tau)\leq C\frac{|\xi|^{2}}
{\langle\sigma\rangle^{-b^{\prime}}\prod\limits_{j=1}^{2}\langle\sigma_{j}\rangle^{b}}
\leq C\frac{|\xi_{1}^{4}-\xi_{2}^{4}|^{\frac{1}{2}}}{\prod\limits_{j=1}^{2}\langle\sigma_{j}\rangle^{b}}.\label{3.023}
\end{eqnarray}

This case can be proved similarly to Case (2) of Lemma 3.1.

\noindent
(3)When$(\xi_{1},\xi,\tau_{1},\tau)\in A_{3}$, which yield $|\xi_{1}|\sim|\xi|,|\xi_{1}|\geq4a$, $|\xi_{2}|>a$,
 then, we consider (\ref{2.04})-(\ref{2.06}), respectively.\\
When (\ref{2.04}) is valid, then, we have
\begin{eqnarray}
&&K_{2}(\xi_{1},\tau_{1},\xi,\tau)\leq \frac{|\xi|^{1+s_{1}-s_{2}+4b^{\prime}}|\xi_{2}|^{-s_{2}+b^{\prime}}}
{\prod\limits_{j=1}^{2}\langle\sigma_{j}\rangle^{b}}\nonumber\\
&&
\leq C\frac{|\xi_{1}|^{\frac{3}{4}}}{\prod\limits_{j=1}^{2}\langle\sigma_{j}\rangle^{b}}\leq C
\frac{|\xi_{1}^{4}-\xi_{2}^{4}|^{\frac{1}{2}}}{\prod\limits_{j=1}^{2}\langle\sigma_{j}\rangle^{b}}.\label{3.024}
\end{eqnarray}
This case can be proved similarly to    Case (2) of Lemma 3.1.\\
When (\ref{2.05}) is valid, which yield $\langle\sigma_{1}\rangle^{-b}\langle\sigma\rangle^{b^{\prime}}\leq
\langle\sigma_{1}\rangle^{b^{\prime}}\langle\sigma\rangle^{-b}$, then, we have
\begin{eqnarray}
&&K_{2}(\xi_{1},\tau_{1},\xi,\tau)\leq \frac{|\xi|^{1+s_{1}-s_{2}+4b^{\prime}}|\xi_{2}|^{-s_{2}+b^{\prime}}}
{\langle\sigma_{2}\rangle^{b}\langle\sigma\rangle^{b}}\nonumber\\
&&\leq C\frac{|\xi_{1}|^{9\epsilon}}{\langle\sigma_{2}\rangle^{b}\langle\sigma\rangle^{b}}\leq
C\frac{|\xi_{1}|^{\frac{3}{4}}}{\langle\sigma_{2}\rangle^{b}\langle\sigma\rangle^{b}}\leq C
\frac{|\xi^{4}-\xi_{2}^{4}|^{\frac{1}{2}}}{\langle\sigma_{2}\rangle^{b}\langle\sigma\rangle^{b}}.\label{3.025}
\end{eqnarray}
This case can be proved similarly to Case (2)  of Lemma 3.1.\\
When (\ref{2.06}) is valid, which yield $\langle\sigma_{2}\rangle^{-b}\langle\sigma\rangle^{b^{\prime}}
\leq \langle\sigma_{2}\rangle^{b^{\prime}}\langle\sigma\rangle^{-b}$, then we have
\begin{eqnarray}
&&K_{2}(\xi_{1},\tau_{1},\xi,\tau)\leq \frac{|\xi|^{1+s_{1}-s_{2}+4b^{\prime}}|\xi_{2}|^{-s_{2}+b^{\prime}}}
{\langle\sigma_{1}\rangle^{b}\langle\sigma\rangle^{b}}\nonumber\\
&&\leq C\frac{|\xi_{1}|^{9\epsilon}}{\langle\sigma_{1}\rangle^{b}\langle\sigma\rangle^{b}}\leq C
\frac{|\xi_{1}|^{\frac{3}{4}}}{\langle\sigma_{1}\rangle^{b}\langle\sigma\rangle^{b}}.\label{3.026}
\end{eqnarray}
This case can be proved similarly to  (2.6) of   Case (3) in Lemma 3.1.\\

\noindent (4)When$(\xi_{1},\xi,\tau_{1},\tau)\in A_{4}$, which yield $|\xi_{1}|\sim|\xi_{2}|\sim|\xi|,|\xi_{1}|\geq4a>2a$,
 then, we consider (\ref{2.04})-(\ref{2.06}), respectively.\\
When (\ref{2.04}) is valid, we have
\begin{eqnarray}
&&K_{2}(\xi_{1},\tau_{1},\xi,\tau)\leq C\frac{|\xi_{1}|^{1+s_{1}-2s_{2}+5b^{\prime}}}
{\prod\limits_{j=1}^{2}\langle\sigma_{j}\rangle^{b}}\nonumber\\
&&\leq C\frac{|\xi_{1}|^{9\epsilon}}{\prod\limits_{j=1}^{2}\langle\sigma_{j}\rangle^{b}}\leq C
\frac{|\xi_{1}|^{\frac{3}{4}}}{\prod\limits_{j=1}^{2}\langle\sigma_{j}\rangle^{b}}\leq C
\frac{|\xi_{1}|^{\frac{3}{8}}|\xi_{2}|^{\frac{3}{8}}}{\prod\limits_{j=1}^{2}\langle\sigma_{j}\rangle^{b}}.\label{3.027}
\end{eqnarray}
This case can be proved similarly to (2.6)  of  Case 4.\\

\noindent When (\ref{2.05}) is valid, which yield $\langle\sigma_{1}\rangle^{-b}
\langle\sigma\rangle^{b^{\prime}}\leq \langle\sigma_{1}\rangle^{b^{\prime}}\langle\sigma\rangle^{-b}$, then, we have
\begin{eqnarray}
&&\hspace{-1cm}K_{2}(\xi_{1},\tau_{1},\xi,\tau)\leq C\frac{|\xi_{1}|^{1+s_{1}-2s_{2}+5b^{\prime}}}
{\langle\sigma_{2}\rangle^{b}\langle\sigma\rangle^{b}}
\leq C\frac{|\xi_{1}|^{9\epsilon}}{\langle\sigma_{2}\rangle^{b}\langle\sigma\rangle^{b}}\leq C
\frac{|\xi_{1}|^{\frac{3}{4}}}{\langle\sigma_{2}\rangle^{b}\langle\sigma\rangle^{b}}\leq C
\frac{|\xi_{2}|^{\frac{3}{8}}|\xi|^{\frac{3}{8}}}{\langle\sigma_{2}\rangle^{b}\langle\sigma\rangle^{b}}.\label{3.028}
\end{eqnarray}
This case can be proved similarly to (2.6)  of  Case 4.\\
When (\ref{2.06}) is valid, which yield $\langle\sigma_{2}\rangle^{-b}\langle\sigma\rangle^{b^{\prime}}
\leq \langle\sigma_{2}\rangle^{b^{\prime}}\langle\sigma\rangle^{-b}$,  we have
\begin{eqnarray}
&&K_{2}(\xi_{1},\tau_{1},\xi,\tau)\leq C\frac{|\xi_{1}|^{1+s_{1}-2s_{2}+5b^{\prime}}}
{\langle\sigma_{1}\rangle^{b}\langle\sigma\rangle^{b}}
\leq C\frac{|\xi_{1}|^{9\epsilon}}{\langle\sigma_{1}\rangle^{b}\langle\sigma\rangle^{b}}\leq C
\frac{|\xi_{1}|^{\frac{3}{8}}|\xi|^{\frac{3}{8}}}{\langle\sigma_{1}\rangle^{b}\langle\sigma\rangle^{b}}.\label{3.029}
\end{eqnarray}
This case can be proved similarly to (2.6)  of  Case 4.\\
(5)When$(\xi_{1},\xi,\tau_{1},\tau)\in A_{5}$, which yield $|\xi_{1}|\sim|\xi_{2}|\sim|\xi|,|\xi_{1}|\geq 4a$.
This case can be proved similarly to Case (4).\\
(6)When$(\xi_{1},\xi,\tau_{1},\tau)\in A_{6}$, which yield $|\xi_{1}|\sim|\xi_{2}|,|\xi_{1}|\geq4a,4|\xi|<|\xi_{2}|$,
 then, we consider (\ref{2.04})-(\ref{2.06}), respectively.\\
When (\ref{2.04}) is valid, we have
\begin{eqnarray}
&&K_{2}(\xi_{1},\tau_{1},\xi,\tau)\leq C\frac{|\xi|^{1+b^{\prime}}\langle\xi\rangle^{s_{1}}|\xi_{1}|^{-2s_{2}+4b^{\prime}}}
{\prod\limits_{j=1}^{2}\langle\sigma_{j}\rangle^{b}}\nonumber\\
&&\leq C\frac{|\xi|^{\frac{1}{2}}|\xi_{1}|^{1+s_{1}-2s_{2}+4b^{\prime}+2\epsilon}}{\prod\limits_{j=1}^{2}\langle\sigma_{j}\rangle^{b}}\leq C\frac{|\xi_{1}^{4}-\xi_{2}^{4}|^{\frac{1}{2}}}{\prod\limits_{j=1}^{2}\langle\sigma_{j}\rangle^{b}}.\label{3.030}
\end{eqnarray}
This case can be proved similarly to Case 2 of Lemma 3.1.\\
When (\ref{2.05}) is valid, which yield $\langle\sigma_{1}\rangle^{-b}\langle\sigma\rangle^{b^{\prime}}\leq
 \langle\sigma_{1}\rangle^{b^{\prime}}\langle\sigma\rangle^{-b}$, then, we have
\begin{eqnarray}
&&K_{2}(\xi_{1},\tau_{1},\xi,\tau)\leq C\frac{|\xi|^{1+b^{\prime}}\langle\xi\rangle^{s_{1}}|\xi_{1}|^{-2s_{2}+4b^{\prime}}}
{\langle\sigma_{2}\rangle^{b}\langle\sigma\rangle^{b}}\nonumber\\
&&\leq C\frac{|\xi|^{\frac{1}{2}}|\xi_{1}|^{1+s_{1}-2s_{2}+4b^{\prime}+2\epsilon}}{\langle\sigma_{2}\rangle^{b}\langle\sigma\rangle^{b}}\leq C\frac{|\xi_{2}^{4}-\xi^{4}|^{\frac{1}{2}}}{\langle\sigma_{2}\rangle^{b}\langle\sigma\rangle^{b}}.\label{3.031}
\end{eqnarray}
This case can be proved similarly to Case 2 of Lemma 3.1.\\
When (\ref{2.06}) is valid, which yield $\langle\sigma_{2}\rangle^{-b}\langle\sigma\rangle^{b^{\prime}}\leq
\langle\sigma_{2}\rangle^{b^{\prime}}\langle\sigma\rangle^{-b}$, then, we have
\begin{eqnarray}
&&K_{2}(\xi_{1},\tau_{1},\xi,\tau)\leq C\frac{|\xi|^{1+b^{\prime}}\langle\xi\rangle^{s_{1}}|\xi_{1}|^{-2s_{2}+4b^{\prime}}}
{\langle\sigma_{1}\rangle^{b}\langle\sigma\rangle^{b}}\nonumber\\
&&\leq C\frac{|\xi|^{1+b^{\prime}}\langle\xi\rangle^{s_{1}}|\xi_{1}|^{-2s_{2}+4b^{\prime}}}{\langle\sigma_{1}\rangle^{b}\langle\sigma\rangle^{b}}\leq C\frac{|\xi|^{\frac{1}{2}}|\xi_{1}|^{1+s_{1}-2s_{2}+4b^{\prime}+2\epsilon}}{\langle\sigma_{1}\rangle^{b}\langle\sigma\rangle^{b}}\leq C\frac{|\xi_{1}^{4}-\xi^{4}|^{\frac{1}{2}}}{\langle\sigma_{1}\rangle^{b}\langle\sigma\rangle^{b}}.\label{3.032}
\end{eqnarray}
This case can be proved similarly to Case 2 of Lemma 3.1.

This completes the proof of Lemma 3.2.

 To prove Theorem 1.3, it suffices to prove Lemma 3.3.

\begin{Lemma}\label{lem3.3}
For $s\leq -\frac{1}{2}$, $b^{\prime}=-\frac{1}{2}+2\epsilon$ and  $b=\frac{1}{2}+\frac{\epsilon}{2}$. Then, we have that
\begin{eqnarray}
&&\left\|\partial_{x}\left(u_{1}u_{2}\right)\right\|_{X_{\frac{1}{2}+2\epsilon,b^{\prime}}}
\leq C\left\|u_{1}\right\|_{X_{s,b}}\left\|u_{2}\right\|_{X_{s,b}}\label{3.033}
\end{eqnarray}
fails.
\end{Lemma}
\noindent {\bf Proof.} We define
\begin{eqnarray*}
&&A=\left\{(\xi,\tau)\in \R^{2}: |\xi-N|\leq N^{-\frac{3}{2}}, \left|\tau-(5\alpha N^{4}+3\beta N^{2})\xi+4\alpha N^{5}+2\beta N^{3}\right|\leq \frac{1}{2}\right\}
,\\
&&B=\left\{(\xi,\tau)\in \R^{2}: |\xi-2N^{-\frac{3}{2}}|\leq N^{-\frac{3}{2}}, \left|\tau-(5\alpha N^{4}+3\beta N^{2})\xi\right|\leq \frac{1}{2}\right\}
,\\
&&R=\left\{(\xi,\tau)\in \R^{2}: |\xi-N|\leq \frac{N^{-\frac{3}{2}}}{4}, \left|\tau-(5\alpha N^{4}+3\beta N^{2})\xi+4\alpha N^{5}+2\beta N^{3}\right|\leq \frac{1}{2}\right\}
,\\
&&f(\xi,\tau)=\chi_{A}(\xi,\tau), g(\xi,\tau)=\chi_{B}(\xi,\tau).
\end{eqnarray*}
Here, $A, B, R$ are defined as in example 2 of \cite{K2011}.
Then, we have
\begin{eqnarray}
f*g(\xi,\tau)\geq CN^{-\frac{3}{2}}\chi_{R}(\xi,\tau).\label{3.034}
\end{eqnarray}
Combining (\ref{3.033}) with (\ref{3.034}), we have
\begin{eqnarray}
N^{\frac{3}{2}+2\epsilon}N^{-\frac{3}{2}}N^{-\frac{3}{4}}\leq CN^{\frac{5b}{2}}N^{-\frac{3}{4}}N^{s}N^{-\frac{3}{4}},\label{3.035}
\end{eqnarray}
which is  equivalent to $s\geq -\frac{1}{2}+\frac{3\epsilon}{4}.$ This contradicts with $s\leq -\frac{1}{2}.$

This ends the proof of Lemma 3.3.

\bigskip
\noindent{\large\bf 4. Proof of Theorem 1.4}
\setcounter{equation}{0}
\setcounter{Theorem}{0}

\setcounter{Lemma}{0}

\setcounter{section}{4}

\begin{Lemma}\label{lem4.1} Let $u_{0}\in H^{s}(\R)(s\geq\frac{1}{4})$. Then, we have
\begin{eqnarray}
&&\lim\limits_{N\rightarrow\infty}\|u-u_{N}\|_{L_{x}^{4}L_{t}^{\infty}}=0.\label{4.01}
\end{eqnarray}
\end{Lemma}
\noindent{\bf Proof.} Inspired by the idea of  the proof process of  Theorem 1.1 of \cite{CLS}.  We firstly prove that the
Cauchy problem for (\ref{1.01}) is locally well-posed in $H^{s}(\R)(s>-\frac{7}{4}).$
We define
\begin{eqnarray}
&&\Phi(u)=\eta(t)U(t)f+\eta\left(\frac{t}{T}\right)\int_{0}^{t}U(t-\tau)\partial_{x}(u^{2})d\tau,\label{4.02}\\
&&B=\left\{u\in X_{s,b}:\|u\|_{X_{s,b}}\leq 2C\|u_{0}\|_{H^{s}}\right\}.\label{4.03}
\end{eqnarray}
By using (\ref{4.01})-(\ref{4.02}) and Lemma 3.1,  for  $T\leq  \left(\frac{1}{4C^{2}\|f\|_{H^{s}}}\right)^{\frac{2}{3\epsilon}}$,     we have
\begin{eqnarray}
&&\left\|\Phi(u)\right\|_{X_{s,b}}\leq \left\|\eta(t)U(t)f\right\|_{X_{s,b}}+
\left\|\eta\left(\frac{t}{T}\right)\int_{0}^{t}U(t-\tau)\partial_{x}(u^{2})d\tau\right\|_{X_{s,b}}\nonumber\\
&& \leq C\|f\|_{H^{s}}+CT^{\frac{3\epsilon}{2}}\left\|\partial_{x}(u^{2})\right\|_{X_{s,b_{1}}}\nonumber\\
&&\leq C\|f\|_{H^{s}}+CT^{\frac{3\epsilon}{2}}\|u\|_{X_{s,b}}^{2}\nonumber\\
&&\leq C\|f\|_{H^{s}}+CT^{\frac{3\epsilon}{2}}(2C\|f\|_{H^{s}})^{2}\leq 2C\|f\|_{H^{s}}\label{4.04}
\end{eqnarray}
and
\begin{eqnarray}
&&\left\|\Phi(u)-\Phi(v)\right\|_{X_{s,b}}\leq
\left\|\eta\left(\frac{t}{T}\right)\int_{0}^{t}U(t-\tau)\partial_{x}(u^{2}-v^{2})d\tau\right\|_{X_{s,b}}\nonumber\\
&&\leq CT^{\frac{3\epsilon}{2}}\|u-v\|_{X_{s,b}}\left[\|u\|_{X_{s,b}}+\|v\|_{X_{s,b}}\right]\leq  2C^{2}T^{\frac{3\epsilon}{2}}\|f\|_{H^{s}}\|u-v\|_{X_{s,b}}\nonumber\\
&&\leq \frac{1}{2}\|u-v\|_{X_{s,b}}\label{4.05}.
\end{eqnarray}
Thus, $\Phi$ is a contraction mapping from $B$ to $B$. Consequently, $\Phi$ has a fixed point.
 That is $\Phi(u)=u.$ From (\ref{4.05}), we have
\begin{eqnarray}
\left\|u-v\right\|_{X_{s,b}}\leq
\frac{1}{2}\|u-v\|_{X_{s,b}}\label{4.06}.
\end{eqnarray}
We
can assume that $u_{N}$ is the solution to the truncated Kawahara equation
\begin{eqnarray}
\partial_{t}u_{N}+\alpha\partial_{x}^{5}u_{N}+\beta\partial_{x}^{3}u_{N}+\partial_{x}\left((u_{N})^{2}\right)=0,\label{4.07}
\end{eqnarray}
with the initial data $P_{N}u_{0}$. We can see that $u_{N}$ is smooth, as the initial
data $P_{N}u_{0}$ is smooth with $u_{0}\in H^{s}(\R)(s\geq\frac{1}{4})$.
Let $u:=u_{\infty}$ be  the solution to the Kawahara equation with initial
  data $u_{0}=P_{\infty}u_{0}$.   We define
\begin{eqnarray}
\Phi(u_{N}):=\eta(t)U(t)P_{N}u_{0}-\eta\left(\frac{t}{T}\right)\int_{0}^{t}U(t-t^{\prime})P_{N}\partial_{x}((u_{N})^{2})dt^{\prime}.\label{4.08}
\end{eqnarray}
Obviously, by using a  proof similar to above, $\Phi$
is a contraction mapping on the ball $\left\{u_{N}:\|u_{N}\|_{X_{s,\frac{1}{2}+\epsilon}}\leq 2C\|u_{0}\|_{H^{s}(\SR)}\right\}$.
 By using Lemma 2.5,  we have
\begin{eqnarray}
&&\left\|\sup\limits_{t\in[0,T]}|u-u_{N}|\right\|_{L_{x}^{4}}\leq C\|u-u_{N}\|_{X_{s,\frac{1}{2}+\epsilon}}(s\geq\frac{1}{4}).\label{4.09}
\end{eqnarray}
 For $t\in [-T,T]$, we have
\begin{eqnarray}
u-u_{N}=\eta(t)U(t)P^{N}u_{0}(x)-\eta\left(\frac{t}{T}\right)\int_{0}^{t}U(t-t^{'})(\partial_{x}(u^{2})-P_{N}\partial_{x}((u_{N})^{2})dt^{\prime}.\label{4.010}
\end{eqnarray}
Then, by using Lemma 2.1, we have
\begin{eqnarray}
&&\|u-u_{N}\|_{X_{s,\frac{1}{2}+\epsilon}}\leq C\|P^{N}u_{0}\|_{H^{s}(\SR)}+CT^{\frac{3\epsilon}{2}}
\|\partial_{x}(u^{2})-P_{N}\partial_{x}(u_{N}^{2})\|_{X_{s,-\frac{1}{2}+2\epsilon}}.\label{4.011}
\end{eqnarray}
Since
\begin{eqnarray}
&&\partial_{x}(u^{2})-P_{N}\partial_{x}(u_{N}^{2})
=P_{N}\left(\partial_{x}(u^{2})-\partial_{x}(u_{N}^{2})\right)+P^{N}(\partial_{x}(u^{2})),\label{4.012}
\end{eqnarray}
 we have
\begin{eqnarray}
&&\|u-u_{N}\|_{X_{s,\frac{1}{2}+\epsilon}}\leq C\|P^{N}u_{0}\|_{H^{s}(\SR)}+CT^{\frac{3\epsilon}{2}}
\|P_{N}\left(\partial_{x}(u^{2})-\partial_{x}(u_{N}^{2})\right)\|_{X_{s,-\frac{1}{2}+2\epsilon}}\nonumber\\
&&+CT^{\epsilon}\|P^{N}(\partial_{x}(u^{2}))\|_{X_{s,-\frac{1}{2}+2\epsilon}}.\label{4.013}
\end{eqnarray}
By using Lemmas 2.1, 3.2 and (\ref{4.07}),      for  $T\leq  \left(\frac{1}{4C^{2}\|f\|_{H^{s}}}\right)^{\frac{1}{\epsilon}}$,           we have
\begin{eqnarray}
&&T^{\frac{3\epsilon}{2}}\|P_{N}\left(\partial_{x}(u^{2})-\partial_{x}(u_{N}^{2})\right)\|_{X_{s,-\frac{1}{2}+2\epsilon}}\nonumber\\&&\leq  CT^{\frac{3\epsilon}{2}}\|P_{N}\left(\partial_{x}[(u+u_{N})(u-u_{N})]\right)\|_{X_{s,-\frac{1}{2}+2\epsilon}}\nonumber\\
&&\leq CT^{\frac{3\epsilon}{2}}\|u+u_{N}\|_{X_{s,\frac{1}{2}+\epsilon}}\|u-u_{N}\|_{X_{s,\frac{1}{2}+\epsilon}}\leq CT^{\frac{3\epsilon}{2}}(\|u\|_{X_{s,\frac{1}{2}+\epsilon}}+\|u_{N}\|_{X_{s,\frac{1}{2}+\epsilon}})\|u-u_{N}\|_{X_{s,\frac{1}{2}+\epsilon}}\nonumber\\
&&\leq 2C T^{\frac{3\epsilon}{2}}\|u_{0}\|_{H^{s}(\SR)} \|u-u_{N}\|_{X_{s,\frac{1}{2}+\epsilon}}\leq \frac{1}{2} \|u-u_{N}\|_{X_{s,\frac{1}{2}+\epsilon}},\label{4.014}
\end{eqnarray}
inserting (\ref{4.014}) into (\ref{4.013}), we have
\begin{eqnarray}
&&\|u-u_{N}\|_{X_{s,\frac{1}{2}+\epsilon}}\nonumber\\&&\leq C\|P^{N}u_{0}\|_{H^{s}(\SR)}+\|P^{N}(\partial_{x}(u^{2}))
\|_{X_{s,-\frac{1}{2}+2\epsilon}}+\frac{1}{2} \|u-u_{N}\|_{X_{s,\frac{1}{2}+\epsilon}}.\label{4.015}
\end{eqnarray}
From (\ref{4.015}), we have
\begin{eqnarray}
&&\|u-u_{N}\|_{X_{s,\frac{1}{2}+\epsilon}}\leq 2C\|P^{N}u_{0}\|_{H^{s}(R)}+
\|P^{N}(\partial_{x}(u^{2}))\|_{X_{s,-\frac{1}{2}+2\epsilon}}.\label{4.016}
\end{eqnarray}
From Lemma 3.1, we have
\begin{eqnarray}
&&\|\partial_{x}(u^{2})\|_{X_{s,-\frac{1}{2}+2\epsilon}}\leq C\|u\|_{X_{s,\frac{1}{2}+\epsilon}}^{2}
\leq 4C^{3}\|u_{0}\|_{H^{s}(\SR)}^{2}<\infty.\label{4.017}
\end{eqnarray}
From  (\ref{4.017}), we have
\begin{eqnarray}
\|P^{N}(\partial_{x}(u^{2}))\|_{X_{s,-\frac{1}{2}+2\epsilon}}\rightarrow 0\label{4.018}
\end{eqnarray}
as $N\rightarrow \infty$.
Since $u_{0}\in H^{s}(\R)(s\geq \frac{1}{4})$, we have
\begin{eqnarray}
\left\|P^{N}u_{0}\right\|_{H^{s}(\SR)}\rightarrow 0\label{4.019}
\end{eqnarray}
as $N\rightarrow \infty$.

Inserting (\ref{4.018}),  (\ref{4.019}) into (\ref{4.016}) yields (\ref{4.01}).

This completes the proof of Lemma 4.1.

To obtain Theorem 1.4, it suffices to prove Lemma 4.2.
\begin{Lemma}\label{lem4.2} Let $u_{0}\in H^{s}(\R)(s\geq\frac{1}{4})$.
Then, we have $u(x,t)\rightarrow u_{0}(x)$ as $t\rightarrow0$ for almost everywhere $x\in \R$.
\end{Lemma}
\noindent{\bf  Proof.}  Inspired by the idea of  the proof of  Proposition 3.3  of \cite{CLS}, we present the proof of Lemma 4.2.
Since $u_{N}$ is smooth,  for all $x\in \R$,  we have
\begin{eqnarray}
&&\lim\limits_{t\rightarrow0}u_{N}(x,t)=P_{N}u_{0}(x).\label{4.020}
\end{eqnarray}
By using the triangle inequality, we have
\begin{eqnarray}
&&|u-u_{0}|\leq |u-u_{N}|+|u_{N}-P_{N}u_{0}|+|u_{0}-P_{N}u_{0}|.\label{4.021}
\end{eqnarray}
By using (\ref{4.021}), we have
\begin{eqnarray}
&&\lim_{t\rightarrow0}\sup|u-u_{0}|\leq\lim_{t\rightarrow0}\sup|u-u_{N}|+|P^{N}u_{0}|.\label{4.022}
\end{eqnarray}
For $\lambda>0$,  by using the Chebyshev inequality and Sobolev embedding as well as (\ref{4.022}), we have
\begin{eqnarray}
&&|\{x\in \R:\lim\limits_{t\rightarrow0}\sup|u-u_{0}|>\lambda\}|\leq |\{x\in \R:\lim\limits_{t\rightarrow0}\sup|u-u_{N}|>\frac{\lambda}{2}\}|\nonumber\\
&&+|\{x\in \R:|P_{N}u_{0}|>\frac{\lambda}{2}\}|\leq C\lambda^{-4}\left\|u-u_{N}\right\|_{L_{x}^{4}L_{t}^{\infty}}^{4}+\lambda^{-2}\|P^{N}u_{0}\|_{L^{2}}\nonumber\\
&&\leq C\lambda^{-4}\left\|u-u_{N}\right\|_{L_{x}^{4}L_{t}^{\infty}}^{4}+\lambda^{-2}\|P^{N}u_{0}\|_{H^{s}}.\label{4.023}
\end{eqnarray}
Since $u_{0}\in H^{s}(\R)(s\geq\frac{1}{4})$, we have
 \begin{eqnarray}
        \|P^{N}u_{0}\|_{H^{s}}\rightarrow0    \label{4.024}
 \end{eqnarray}
 as $N\rightarrow\infty$.  By using Lemma 4.1 and (\ref{4.024}),   we have
\begin{eqnarray*}
&&|\{x\in \R:\lim\limits_{t\rightarrow0}\sup|u-u_{0}|>\lambda\}|=0.
\end{eqnarray*}

This completes the proof of Lemma 4.2.

\bigskip
\noindent{\large\bf 5. Proof of Theorem 1.5}
\setcounter{equation}{0}
\setcounter{Theorem}{0}

\setcounter{Lemma}{0}

\setcounter{section}{5}
In this section, we use Lemmas 2.1,  3.1 to prove Theorem 1.5.

\noindent{\bf  Proof of Theorem 1.5.} Inspired by page 7 of \cite{Compaan}, we present the proof of Theorem 1.5.
Obviously,
\begin{eqnarray}
&&u(x,t)-U(t)u_{0}(x)=\eta\left(\frac{t}{T}\right)\int_{0}^{t}U(t-t^{\prime})\partial_{x}(u^{2})dt^{\prime}.\label{5.01}
\end{eqnarray}
Let $s_{1}=\frac{1}{2}+2\epsilon$, $b=\frac{1}{2}+\epsilon$, $s_{2}\geq-\frac{1}{2}+\epsilon$. By using Lemmas 2.1,  3.1  and  (\ref{5.01}), we have
\begin{eqnarray}
&&\|u-U(t)u_{0}\|_{X_{s_{1},b}}= \left\|\eta\left(\frac{t}{T}\right)\int_{0}^{t}U(t-t^{\prime})(\partial_{x}(u^{2}))dt^{\prime}\right\|_{X_{s_{1},b}}\nonumber\\
&&\leq CT^{\frac{3\epsilon}{2}}\|\partial_{x}(u^{2})\|_{X_{s_{1},-\frac{1}{2}+2\epsilon}}\leq C\|u\|_{X_{s_{2},b}}^{2}\leq 2C^{3}\|u_{0}\|_{H^{s_{2}}(\SR)}^{2}<\infty.\label{5.02}
\end{eqnarray}
Thus, from Lemma 2.7 and $H^{s_{1}}(\R)\hookrightarrow C(\R)(s_{1}=\frac{1}{2}+2\epsilon)$,  we have
 $u-U(t)u_{0}\in X_{s_{1},\>b}\hookrightarrow C([-T,T];H^{s_{1}}(\R))\hookrightarrow C([-T,T];C(\R))$.
 Thus,   we have
\begin{eqnarray}
&&\lim\limits_{t\rightarrow0}\sup\limits_{x\in\SR}|u(x,t)-U(t)u_{0}|=0.\label{5.03}
\end{eqnarray}

This ends the proof of Theorem 1.5.

\bigskip

\leftline{\large \bf Acknowledgments}
 This work is supported by the education department of Henan Province under
grant number 21A110014.

\end{document}